\documentclass[12pt]{article}
\usepackage{amsmath,amsthm,amsfonts}
\input xy
\xyoption{all}
\newtheorem{theorem}{Theorem}[section]
\newtheorem{lemma}[theorem]{Lemma}

\newtheorem{corollary}[theorem]{Corollary}
\newtheorem{claim}{Claim}

\newtheorem*{ex}{Example}
\renewcommand{\P}{\mathbb{P}}
\newcommand{\C}{\mathbb{C}}
\begin{document}
\title{On singular varieties having an extremal secant line}
\author{Marie-Am\' elie Bertin\\
12 rue des Bretons,\\
94700 Maisons-Alfort\\
FRANCE \\
\texttt{marie.bertin@parisfree.com}\\
\texttt{mabertin@math.jussieu.fr}
}

\date{\today}
\maketitle

\section{Introduction}
Let $X$ be an irreducible non degenerate variety of degree $d$ and codimension $e$ in ${\P}^r$. 
Let us recall that, if $X$ is smooth, the regularity conjecture
foresees that the Castelnuovo-Mumford regularity of $X$ is less or equal to $d-e+1$. This conjecture has been proved to hold 
for curves without the smoothness assumption by L.Gruson, R.Lazarsfeld and C.Peskine \cite{GLP} and 
for smooth surfaces by H.Pinkham \cite{P} (for the normality part of the conjecture) and R.Lazarsfeld \cite{L}. 
In a previous article \cite{B}, we've extended Gruson, Lazarsfeld and Peskine's result to smooth scrolls over curves.
 In the curve case, Gruson, Lazarsfeld and Peskine give moreover a complete classification of curves of extremal regularity; 
they all posses an extremal secant line but the elliptic normal curves and the singular rational curves of degree $d$ in ${\P}^{d-1}$.
In case the regularity conjecture holds, this leads us to believe that smooth varieties of extremal regularity of fixed dimension
have an extremal secant line but a finite list of varieties. Varieties with an extremal secant line are moreover $(d-e)$-irregular,
so they form a good test of the conjecture. This has lead us to classify the smooth ones and compute their regularity in a previous article \cite{B}.
We gave a complete classification of the smooth varieties with an extremal secant line if $e\geq 2$;
they are rational scrolls with an extremal secant line, the Veronese surface in ${\P}^5$ or its smooth projection to ${\P}^4$. 
This classification can be seen as a generalization of del Pezzo and Bertini 's classification
of varieties of minimal degree, since such varieties correspond to varieties for which an extremal secant line is a bisecant line.
Our classification results have been extended to a classification of smooth varieties
having an extremal or next to extremal secant linear space $L$ by S.Kwak \cite{K}, 
under the extra assumption the extremal secant linear space $L$ cuts $X$ along a curvilinear 
$0$-dimensional scheme, i.e. $dim(L\cap \mathbb{T}_{p}(X))\leq 1$ for all points $p\in L\cap X$. 
Independently, A. Noma \cite{N} has recently  found a generalization of Kwak's result:

\begin{theorem}[B. ($x=0$, $k=1$), Kwak($x\leq 1$), Noma] 
Let $X$ be a smooth, complex, irreducible, non degenerate variety  of degree $d$ and  codimension $e$ of ${\P}^r$.
Let $L$ be a $k$-dimensional linear space intersecting $X$ in
a $0$-dimensional scheme of length $m$, such that $L\cap X$ is curvilinear.
 Assume that both $m\leq d-e+k-x$ and $e-k\geq x+1$.
Then the sectional genus $\pi (X)$ satisfies $\pi (X)\leq x$.
\end{theorem}

Let us point out that Kwak's result (and ours) is more precise for it gives a classification result.
In this note we'll see that this bound still holds for $x=0$ and $k=1$ without the smoothness assumption.
In view of a classification theorem of varieties having an extremal secant line in the singular case, 
our previous paper \cite{B} also stated a key result, that is unfortunately false in the singular case. 
In this note we correct it and deduce a classification of singular varieties of codimension $e\geq 2$ 
having an extremal secant line.

The first of these two key results in the singular case is: an extremal secant line to $X$ meets $X$  along smooth points. 
This is wrong; A. Noma kindly pointed out to us a very simple counterexample. 

\begin{ex}[Noma's counterexample]
Let $C$ be a rational normal cubic curve and $X:=<q,C>$ a cone over $C$ in ${\P}^4 =<q,<C>>$.
A general line $l$ through $q$ meets $X$ with multiplicity $2$ at $q$ and is an extremal secant line to $X$. 
\end{ex}
In this example one can find nonetheless an extremal secant line meeting $X$ along smooth points of $X$, namely, any bisecant line to $C$.
This is a general phenomenon and we show in this note that we can correct our key lemma as follows:

\begin{theorem} Let $X$ be an irreducible, non degenerate variety of degree $d$ dimension $n\geq 2$ and codimension $e\geq 2$ in ${\P}^r$.
Suppose that $X$ admits an extremal secant line $l$. Then, there exists an extremal secant line $l^{\prime}$ meeting $X$ along smooth points of $X$.
\end{theorem}

We deduce from this result the following classification theorem.

\begin{theorem} 
Let $X$ be as in the previous theorem. Suppose that $e\geq 2$. Then the variety $X$ is either
\begin{enumerate}
\item a cone $<L,V>$ over the Veronese surface $V$ in ${\P}^5$, where $L$ is a linear space of dimension $k\geq -1$

\item a cone $<L,V^{\prime}>$ over $V^{\prime}$, the isomorphic projection of $V$ to ${\P}^4$, where $L$ is a linear space of dimension $k\geq -1$,

\item a cone $<L,X_0>$, where $X_0$ is a smooth rational scroll with an extremal secant line or a smooth rational curve with an extremal secant line, 
and $L$  a linear space of dimension $k\geq -1$.

\end{enumerate}
\end{theorem}

\section{Basic facts and notations}

In this section we set up notations and recall the main results we use in the sequel.
Let $X$ be a $n$-dimensional complex projective non-degenerate irreducible variety of degree $d$ in ${\P}^{r}$. We set $e:=r-n$, the codimension of $X$.
In the rest of this article, the bracket $<\cdot >$ denotes the linear span  of the subvarieties of $\mathbb{P}^r$ listed in it.
 
A $k$-secant $m$-secant plane to $X$ is a $m$-plane $L$ such that $L \cap X$ is a $0$-dimensional scheme of degree at least $k$. We have the following 
classical result bounding $k$:

\begin{theorem}[linear section theorem] \label{th:sec} Let $X$ be as above. 
Let $\Lambda$ be a linear subspace of ${\P}^r$ of dimension $s\leq e$. Suppose that $\Lambda$ is not contained in $X$, then
it cuts $X$ along a zero dimensional scheme of degree at most $d-e+s$.
\end{theorem}

A proof can be found for instance in \cite{K}.

In particular, a $k$-secant line $l$ to $X$ satisfies $k\leq d-e+1=:\delta$, so that we say that $l$ is an extremal secant line to $X$ if $k=d-e+1$. 
For example, varieties of minimal degree $Z$, i.e. such that $d=e+1$, of degree $d\geq 2$ have an extremal secant line; indeed, any bisecant line 
to $Z$ is an extremal secant line to $Z$.

Recall that a rational normal scroll of dimension $n\geq 2$ is the tautological embedding of some rank $n$ projective bundle over ${\P}^1$,
say $\P (\mathcal{O}_{{\P}^1}(a_1)\oplus \cdots \oplus \mathcal{O}_{{\P}^1}(a_n))$  for some integers $0\leq a_1 \leq \cdots \leq a_n$.
We denote such a scroll by $S(a_1,\cdots ,a_n)$; it is a variety of degree $d=a_1+ \cdots +a_n$, whose linear span is  ${\P}^{d+n-1}$. 
It is therefore a variety of minimal degree. If $d>1$, the variety $S(a_1, \cdots ,a_n)$ is smooth provided that $a_1 >0$ .

In case $n=2$, for all integer $k\geq 0$, we denote by $\mathbb{F}_{k}$ the projective ${\P}^1$-bundle $\P (\mathcal{O}_{{\P}^1}\oplus \mathcal{O}_{{\P}^1}(k))$.
    
We have the following classification of varieties of minimal degree due to del Pezzo and Bertini.

\begin{theorem}[del Pezzo, Bertini] \label{th:min} Let $X$ be an irreducible, non-degenerate variety of codimension $e$ and degree $e+1$ in ${\P}^r$, i.e.
bisecant lines are extremal secant lines to $X$. Then $X$ is a cone $<L,X_0>$, where $L$ is linear space of dimension $k\geq -1$ and
$X_0$ is either 

\begin{enumerate}
\item ($e=1$) a linear space,

\item or  ($e=2$) a smooth quadric hypersurface,

\item or a smooth rational normal scroll,

\item or or the Veronese surface in ${\P}^{5}$.
\end{enumerate} 
\end{theorem}
A modern proof of this can be found in Eisenbud and Harris survey \cite{EH}.

Let $\mathcal{I}_{X|{\P}^r}$ be the ideal sheaf of $X$ in ${\P}^r$. The variety $X$ is said to be $k$-regular, for $k \in \mathbb{Z}$, 
if the following vanishing occur
\[
H^{i}(\mathcal{I}_{X|{\P}^r}(k-i))=0 \, , \qquad \forall i\geq 1
\]

The $k$-regularity property of $X$ implies the $(k+1)$-regularity of $X$, so that one can define the Castelnuovo-Mumford regularity of $X$
as
\[
reg(X):=min\{ k\in \mathbb{Z} | X \,\,\text{is}\,\, k-\text{regular}\}
\]

If $X$ is $k$-regular, the saturated ideal, $I_{X|{\P}^{r}}$, defining $X$ in ${\P}^r$ is generated by polynomials of degree $\leq k$;
thus, if $X$ has a $k$-secant line, $reg(X)\geq k$.

The regularity conjecture (\cite{Ca},\cite{EG}) foresees that:
$reg(X)\leq d-e+1$ (Castelnuovo's bound). 

This conjecture has been proved to hold for curves by Gruson, Lazarsfeld and Peskine \cite{GLP}. In this beautiful paper, they also completely
classify $(d-e)$-irregular curves. This is a key result for the classification of varieties with an extremal secant line.

\begin{theorem}[Gruson-Lazarsfeld-Peskine, 1983] \label{th:GLP} Let $C$ be an irreducible non degenerate curve of degree $d$ in ${\P}^r$.
Suppose that the regularity of $C$ is exactly $d-r+2$, then $C$ is either

\begin{enumerate}
\item a smooth rational curve having an extremal secant line,
\item or an elliptic normal curve 
\item or a singular rational curve of degree $d$ in ${\P}^{d-1}$.

\end{enumerate}
In the last two cases, $C$ has no extremal secant lines, i.e no $3$-secant lines.
\end{theorem}

This has led us to check the regularity conjecture on varieties having an extremal secant line \cite{B}.
We have shown that smooth (not necessarily rational) scrolls satisfy the conjecture, using a similar argument to \cite{GLP}.
In particular, we have deduced that smooth varieties having an extremal secant line satisfy the conjecture.
We show in this article that the same conclusion holds for singular varieties of codimension $e\geq 2$ having an extremal secant line.

To achieve this result, we'll need the following result of ours \cite{B}, which holds provided $X$ has an extremal 
secant line $l$ meeting $X$ at smooth points.

\begin{theorem}\label{th:B00} Let $X$ be an irreducible, non degenerated variety of degree $d$ and codimension $e$ in ${\P}^r$.
Suppose that $X$ has an extremal secant line $l$, such that $l$ meets $X$ along smooth points of $X$.
Then $X$ is the regular projection of a variety of minimal degree $\overline{X}$ of degree $d$ in ${\P}^{d+n-1}$.
\end{theorem}
  
Finally, we recall Bertini's famous irreducibility theorem in the form we shall use. 

\begin{theorem}[Bertini's irreducibility theorem]\label{th:irr} Let $X$ be a $n$-dimensional complex projective subvariety of ${\P}^{r}$
and $\mathcal{L}$ a \underline{non empty} linear system on $X$.
Suppose that $\mathcal{L}$ satisfies the following conditions
\begin{enumerate}
\item there exists a $n$-dimensional linear system $\mathcal{M}$ on ${\P}^r$ such that $\mathcal{L}=\mathcal{M}|_{X}$ 
(if so one can find $\mathcal{M}$ satisfying moreover $dim(\mathcal{M})=dim(\mathcal{L})$),

\item the image of $X$ by the rational map $\Phi_{\mathcal{L}}: X\dashrightarrow {\P}^n$ 
associated to $\mathcal{L}$ is a \underline{non degenerated} subvariety ${\P}^n$ of dimension at least $2$.  
\end{enumerate}
Then, the generic member of $\mathcal{L}$ is an irreducible and non multiple. 
\end{theorem}

This theorem follows from Bertini's theorem on hyperplanes sections of $\Phi_{\mathcal{L}}(X)$ in ${\P}^n$ 
(for a modern proof see for instance \cite{Da} p 249) by duality.
Remark that, if $X$ is non degenerated in ${\P}^r$ and $\Phi_{\Lambda}$ is a rational map on $X$ induced by a linear projection $\pi$, $\mathcal{L}$
satisfies the assumptions of theorem \ref{th:irr} provided that $dim(\pi (X))\geq 2$.
 
\section{The linear system of hyperplanes containing an extremal secant line to $X$}

We need first to establish the following theorem which  determines for which variety $X$ the image of $X$ by the linear system
of hyperplanes passing through $l$ has dimension $<=1$, so that theorem \ref{th:irr}  does not apply.

\begin{theorem}\label{th:1} Let $S$ be a non degenerate, irreducible surface of degree $d$ in ${\P}^{r},r\geq 4$.
Suppose that $l$ is an extremal secant line to $S$, such that the image of $S$ by the linear projection from $l$ is a curve $C$ in ${\P}^{r-2}$.
Then $S$ is a cone $<p,\zeta>$, where $p$ is a point and $\zeta$ is a smooth rational curve of degree $d$ in ${\P}^{r-1}$ 
having an extremal secant line $l^{\prime}$.
\end{theorem}

\begin{proof}
Let $q$ be a generic point on $l$; we denote by $\pi_q$ the projection from $q$.
Let $q_0\in {\P}^{r-1}:=<q_0 , <C>>$ denote the image of $l$ by $\pi_q$.
By dimension count, $\pi_q (S)=<q_0,C>$. Remark also that $S$ lies on the $3$-dimensional cone $Y:=<l,C>$, as a sub-ruled surface,
fibered  in curves of degree $m:=deg(\pi_q )$. We have $deg(C)=d/m$.

\begin{lemma} The projection $\pi_q$ is generically one to one.
\end{lemma}  
\begin{proof}
We wish to show that $m=1$.
Let us desingularize $Y$ and compute the class of the strict transform of $S$ in this desingularization.
Recall that the scroll $Y$ is the image of
\[
\xymatrix{
{\tilde{Y}:=\{ (t,x)\in C\times {\P}^r | x\in <l,t>\}}&{ \subset} &{C\times {\P}^{r}}\ar[dl]_{\phi_1}\ar[dr]^{\phi_2}&\\
&{\P}^1&&{{\P}^r}\\
}
\]
by the projection $\phi_2$ to ${\P}^r$. Let $\phi_1$ denote a generic linear projection of $C\subset {\P}^{r-2}$ onto ${\P}^1$. 
Let $\sigma$ be a fixed isomorphism ${\P}^1\xrightarrow{\sigma} l$.
Consider the graph $\Delta:=\{ (t,x)\in C\times l| \sigma\circ \phi_1 (t)=x\}$;
the projection map $\phi_2$ induces a degree $deg(C)$-map from $\Delta$ onto $l$.  
If $C$ is smooth, $\tilde{Y}\xrightarrow{p} Y$, is a desingularization of $Y$  with exceptional locus
$E:=C\times l$. Let $\overline{C}\xrightarrow{\nu} C$ be the normalization of $C$. 
The proper morphism ${\overline{C}}\times{\P}^r \xrightarrow{(\nu ,id_{{\P}^r})}C\times {\P}^r$, induces a desingularization
$\overline{Y}\xrightarrow{g}\tilde{Y}$ of $\tilde{Y}$, where $\overline{Y}$ denotes the strict transform of $\tilde{Y}$ by $(\nu , id_{{\P}^r})$.
We have the following diagram 
\[
\xymatrix{
{\overline{Y}}& {\subset} &{\overline{C} \times {\P}^{r}} \ar[dl]_{{\overline{\phi_1}}}\ar[dr]^{{\overline{\phi_2}}}&\ar[r]^{(\nu ,id_{{\P}^r})}&&
{C\times{\P}^r}\ar[dl]_{\phi_1}\ar[dr]_{\phi_2}& {\supset}&  {\tilde{Y}}\\
&{\overline{C}}&&{{\P}^r}&C&&{{\P}^r}&\\
}
\]
We have the following equations relating the projection maps:
\begin{equation}\label{Eq:1}
\overline{\phi_2}=\phi_2 \circ (\nu ,id_{{\P}^r})
\end{equation}
and
\begin{equation}\label{Eq:2}
\nu\circ \overline{\phi_1}=\phi_1 \circ (\nu ,id_{{\P}^r})
\end{equation}
It follows that $g$ satisfies $\overline{\phi_2}|_{\overline{Y}} =\phi_2 \circ g$. Since $\overline{\phi_2}|_{\overline{Y}}$ and 
$\phi_2$ are projective, hence proper, the morphism $g$ is proper.
The map $\overline{g}:=g\circ \phi_2$ induces a desingularization of $Y$.  
The Chow ring of $\overline{Y}$ is generated by
$f$, the class of a fiber of the scroll $\overline{Y}\xrightarrow{\overline{\phi}_1} \overline{C}$, 
and $h$ the class of restriction to $\overline{Y}$ of the pull back by $\overline{\phi_2}$ of an hyperplane in ${\P}^r$.
From equation \ref{Eq:2}, we deduce that $f$ is also the class of the pull back by $g$ of a generic fiber of $\phi_1$.
Moreover, equation \ref{Eq:1}, shows that $h$ is also the class of the pull back by $\overline{g}$ of hyperplane sections of $Y$ in ${\P}^r$.
Those generators satisfy the obvious relations: $h^4=0,\, h^3f=0,\, f^2=0$.
Let $\overline{S}$ be the strict transform of $S$ by $\overline{g}$. We have $[\overline{S}]=\alpha h +\beta f$, for some $\alpha,\beta\in\mathbb{Z}$
. Since $g$ is proper and $\phi_2$ is projective, $\overline{g}$ is proper. The projection formula applied to $\overline{g}$, therefore shows
\[
\begin{matrix}
h^3-deg(C)h^2f=0\\
[\overline{S}]\cdot hf=m\\
[\overline{S}]\cdot h^2=d
\end{matrix}
\]
We thus find $\alpha=m$ and $\beta =0$. 
Let us compute the class of  $\overline{\Delta}$, 
the strict transform of $\Delta$ by $g$.
Applying the projection formula to $g$, we get $[\overline{\Delta}]\cdot f=1$, so that $[\overline{\Delta}]=h^2+bhf$ for some $b\in\mathbb{Z}$.
Moreover, applying the projection formula to $\overline{g}$, we get  $[\overline{\Delta}]\cdot h=(deg(C))[l]\cdot h_Y$, 
where $h_Y$ denote the class of an hyperplane section of $Y$ in $Chow({\P}^r)$. 
Thus, $b=deg(C)\mu$, where $\mu$ is the multiplicity at the point $q$, generic on $l$, of the intersection of $l$ with the cone $<q,C>$.
From the relation $\overline{\Delta}\cdot \tilde{S}=d-e+1$, we deduce that
$m(1+\mu )= d-e+1$.  Remark now that $\mu\leq deg(C)-(e-1)+1$. Suppose that l is neither an extremal secant line nor a next to extremal secant line 
to $<q,C>$; we find  $m(deg(C)-(e-1))\geq d-e+1$, so that $m(e-1)\leq e-1$. 
Since $e\geq 2$, it follows that $\delta=1$. Suppose that $l$ is a next to  extremal secant line to $<q,C>$; 
we get $m(e-2)=e-1$, this is clearly impossible.
Finally, if $l$ is an extremal secant line to $<q,C>$, we get $m(e-3)=e-1$. This has clearly no solutions for $e=3$.
If $e\not=3$ we find $\delta=1+\frac{2}{e-3}$, this is also impossible. 
Therefore, we have $m=1$, that is to say $\pi_q$ is birational.
\end{proof}
From $m=1$, we deduce that $S$ is ruled in lines, $\gamma_t$, over $C$. 
\begin{lemma} The intersection $<q,C>\cap X$ has dimension $1$.
Let $\zeta$ be one of the $1$-dimensional irreducible components of $<q,C>\cap S$.
The curve $\zeta$ projects from $q$ onto $C$.
\end{lemma}
\begin{proof}
We have $<q,C>\cap X=\cup_{t\in C} <q,t>\cap \gamma_t$. Therefore, $dim(<q,C>\cap S)=1$.
Let $\zeta$ be a $1$-dimensional irreducible component of $<q,C>\cap S$. 
Suppose, to the contrary, that $\pi_{q_0}(\zeta)$ is a point. 
By dimension count and irreducibility of $\pi_{q_0}(\zeta )$,  we find $\zeta =<q_0 ,\pi_{q_0}(\zeta )>$.
Since $\zeta\subset <q,C>$, $q_0$ lies on $<q,\pi_q (S)>$. Thus, $S\subset<q,q_0,C>=<q,\pi_q (S)>\subset {\P}^{r-1}$. 
This contradicts the non degeneracy of $S$.
\end{proof}
The family of lines $\gamma_t$ can be therefore re-parameterized by $\zeta$. 
Let us denote by $\{p_1 ,\cdots ,p_s \}$ the support of $l\cap S$.
For $i\in \{1,\cdots s\}$, any point $p_i$ of $S\cap l$ must lie on a line $\gamma_t$, for some $t\in C$. 
It follows that $s=1$. Indeed, $S$ has dimension $2$, so the following set is not empty: 
\[
\mathcal{S}:=\{ p\in \{ p_1,\cdots ,p_s \} | \exists \,\,  U \subset X^{\prime\prime}\,\,\textrm{open, dense such that}\,\, \forall t\in U, \,p\in \gamma_t \}
\]
Suppose that $p_1$ belongs to $\mathcal{S}$. By irreducibility assumption, the variety $S$ must coincide with  the cone $<p_1,\zeta>$,
since they intersect along a $2$-dimensional locus. 
The Castelnuovo-Mumford regularity of
$S$, $reg(S)$, satisfies $reg(S)\geq d-e+1$, since $S$ has an extremal secant line. 
Besides, $reg(S)=reg(\zeta)$ for $X$ is a cone.
Therefore, $reg(S)=d-e+1=reg(\zeta)$ and $\zeta$ is a curve of extremal regularity in ${\P}^{r-1}$.
Applying theorem \ref{th:GLP}, we deduce that $\zeta$ is a smooth rational curve having an extremal secant line
, an elliptic normal curve or a singular rational curve of degree $d$ in ${\P}^{d-1}$ ($r=d$).

\begin{lemma} Suppose that $S$ in ${\P}^r$ is a cone $<p,C>$, where $C$ is a curve in ${\P}^{r-1}$. 
Then, there exists an extremal secant line $l$ to $S$ through $p$ if and only if $C$ is a smooth rational curve having an extremal secant line.    
\end{lemma}
\begin{proof}
Let $l$ be a secant line to $S$ through $p$. 
Let $x$ be the intersection of $l$ and $<C>$, $x\not\in C$. We can assume that $p$ has coordinates $(0:\cdots :0:1)$ and $x$ coordinates
$(1:0:\cdots :0)$ in ${\P}^{d}$, so that $I_l:=(x_1 ,\cdots ,x_{d-1})$ is the saturated ideal defining  the line $l$ in the polynomial ring 
$\mathbb{C}[x_0,\cdots ,x_{d}]$. Let $I_C:=(f_1,\cdots ,f_t)$ be the saturated ideal defining $C$ in ${\P}^{d-1}$; the scheme $l\cap S$ has defining ideal 
$(I_C+(x_1 ,\cdots ,x_{d-1}))=(f_1^{0}, \cdots f_{d-1}^{0})+I_l$, where, for $i=1,\cdots , d-1$, the polynomial $f_{i}^{0}$ denotes the leading term
of $f_i$ with respect to the lexicographic monomial ordering in $\mathbb{C}[x_0,\cdots ,x_{d}]$. Therefore, the scheme $l\cap S$ has defining ideal $(f_1^{0},\cdots ,f_{d-1}^{0})$ in $\mathbb{C}[x_0]$. The multiplicity of intersection of $l$ and $S$ at $p$ is therefore $min\{ d_i | f_{i}^{0}\not = 0\}$, 
where $d_i$ denote the degree of the polynomial $f_i$, for $i=1,\cdots ,d-1$.

If $C$ is an elliptic normal curve or a singular curve of degree $d$ in ${\P}^{d-1}$, all the polynomials $f_i$ have degree $2$ 
(see for instance \cite{GLP}: $C$ satisfies property ``$C_2$'' and is non degenerate), so that $l$ is only a bisecant line, while $d-r+2=3$. 
If $C$ is a rational curve having an extremal secant line, we have $d_i\geq d-e+1$, so that $l$ is an extremal secant line to $S$.
\end{proof}

\end{proof}
\begin{corollary}\label{cor:1} Let $X$ be a complex, irreducible non degenerate variety of degree $d$ and dimension $n\geq 2$ in ${\P}^r$.
Assume that $X$ has an extremal secant line $l$. Moreover, suppose that the image of $X$ by the projection $\pi_l$ from  $l$
is a variety $X^{\prime\prime}$ of dimension $n-1$. Then $X$ is a cone $<p,\zeta >$ over a $(n-1)$-dimensional variety $\zeta$ of ${\P}^{r-1}$. 
\end{corollary}

\begin{proof}
Let $X$ in ${\P}^{r}$ be a variety satisfying the assumptions of theorem \ref{th:1}. 
We denote by $p_1,\cdots ,p_s$ the support of $X\cap l$.

\begin{lemma} Let $X$ be a variety satisfying the assumption of  corollary \ref{cor:1}.
By induction on $n\geq 2$, the following properties hold:
\begin{enumerate}
\item $s=1$
\item For $q\in l\setminus \{ p_1 \}$ generic, the projection $\pi_q$ induces a birational map from $X$ onto its image.
\end{enumerate}
\end{lemma}
\begin{proof}
For $n=2$, the result follows from theorem \ref{th:1}.
Suppose that the lemma holds for $(n-1)$-dimensional varieties such that $2\leq n$.
By assumption, the linear system of hyperplanes through $l$ cut on $X$ a linear system satisfying the assumptions of
theorem \ref{th:irr}, hence a generic hyperplane section $h_1$ of $X$ by a hyperplane $H$ through $l$ is an irreducible
variety of dimension $n-1$, and degree $d$ to which $l$ is an extremal secant line. 

By induction hypothesis, we deduce that $s=1$ and that ,for $q\in l\setminus \{ p_1\}$,  $\pi_q$ induces a birational map from $h_1$ onto its image.
Let $t$ be a point of $X^{\prime\prime}\cap H$. This hyperplane $H$ contains the $2$-plane 
$<t,l>$, hence either it induces an hyperplane section of the fiber curve $\gamma_t$ of $X$ or $\gamma_t$ coincides with the line $H\cap <t,l>$.
Let $q\in l\setminus\{ p_1\}$. By induction hypothesis, the projection $\pi_q$ induces a birational map from $h_1:=X\cap H$ onto its image, 
so that the ruled variety $h_1$ is ruled in lines over $\pi_l (h_1)$. The line $H\cap <l,t>$ is one of the lines of the ruling of $h_1$. We deduce that
$\gamma_t =H\cap <l,t>$, so that $\pi_q$ induces a birational map from $X$ onto its image.
\end{proof}
As in the $2$-dimensional case, reparameterizing $X$ by any maximal dimensional component $\zeta$ of $<q,X^{\prime\prime}\cap X$ shows that $X$ is the cone 
$<p_1, \zeta>$.
\end{proof}

\section{Classification of varieties having an extremal secant line}

Going back to the inductive argument of lemma 1.1 in \cite{B}. We can correct it as follows.

\begin{theorem}\label{th:2}
Let $X$ be a complex, irreducible, non-degenerate variety of degree $d$ and codimension $e$ in ${\P}^r$, such that $X$ has an extremal secant line. 
Then, $X$ has an extremal secant line $l$ meeting $X$ along smooth points.
\end{theorem}
 
\begin{proof}
 
Let us proceed by induction on $n$. By theorem \ref{th:GLP}, the result is true for $n=1$.
The result holds for $n=2$ by theorem \ref{th:1}.

Suppose that the result holds for all integers $k$ such that $2\leq k\leq n-1$. 
Let $X$ be  a complex, irreducible, non-degenerate variety of degree $d$ and dimension $n$ 
in ${\P}^r$, such that $X$ has an extremal secant line $l$. 
Consider the linear system of hyperplane sections of $X$ by hyperplanes containing $l$.

Since $n\geq 3$,  we can apply theorem \ref{th:irr} to this linear system. Its general member $h_1$ is thus
irreducible and non multiple. The line $l$ is an extremal secant line to $h_1$. 
Indeed, let $H$ be the hyperplane cutting out $h_1$ on $X$.
By construction, we have the following scheme inclusion: $H\cap X\supseteq l\cap X$. Thus, $l\cap X =l\cap h_1$.

Repeating the same argument for $h_1$, we find that
the generic section $h_{n-2}$ of $X$ by  $(e+2)$-planes containing $l$ is irreducible. The line $l$ is an extremal secant line for the surface $h_{n-2}$.

If the linear system of hyperplanes sections of $h_{n-2}$ passing through $l$ satisfies the assumptions of theorem \ref{th:irr}
, by theorem \ref{th:GLP}, the line $l$ meets a generic hyperplane section of $h_{n-2}$ by $l$ at smooth points. 
Therefore arguing as in \cite{B}, we deduce that $l$ meets $h_{n-2}$ at smooth points. 
Repeating the argument, the line $l$ meets $X$ along smooth points.

If $l$ doesn't meet $X$ along smooth points, the general surface section $h_{n-2}$ of $X$ by  $(e+2)$-planes containing $l$ is a cone $<p,C>$, 
where $C$ is a smooth rational curve of degree $d$ in ${\P}^{r-1}$ having an extremal secant line $l^{\prime}$.
The line $l^{\prime}$ is an extremal secant line for $X$, since $l^{\prime}\cap h_{n-2}\subset l^{\prime}\cap X$ and $length(l^{\prime}\cap h_{n-2})=d-e+1$. 
Repeating the previous argument with $l^{\prime}$ instead of $l$, we deduce that $l^{\prime}$ meets $X$ along smooth points.
\end{proof}

We can now correct proposition 1 in \cite{B}. The argument is as in \cite{B}, the crucial point is, of course, 
the irreducibility of the hyperplane section we use.  
\begin{corollary}\label{cor:2} Let $X$ be an irreducible non degenerate variety of degree $d$ and codimension $e\geq 2$ in ${\P}^{r}$. 
Suppose that $X$ has an extremal secant line $l$ meeting $X$ along smooth points of $X$. 
The image of $X$ by the projection $\pi_l$ from the line $l$
is a variety of minimal degree $X^{\prime\prime}$ in ${\P}^{r-2}$.
\end{corollary}
\begin{proof} Let us proceed by induction on $n:=dim(X)$. 
Notice that, to prove the first assertion, it is enough to show that $dim(X^{\prime\prime})=n$.
Indeed, if it is so, the degree $d^{\prime\prime}$ of $X^{\prime\prime}$ satisfies :
$d^{\prime\prime}\leq d-(d-r+n+1)=(r-2)-n+1$.  
Suppose that $n=1$.Then, the property is trivially true, since $X$ is then a non degenerate non plane curve ($e\geq 2$).
Suppose that $n\geq 2$. By assumption, the linear system of hyperplane sections through $l$ does satisfy Bertini's irreducibility theorem, 
since $X$ meets $l$ along smooth points. Therefore a generic member of the linear system of hyperplanes containing $l$, cuts $X$ 
along an irreducible non multiple variety. This linear system contains a pencil, since $r-2\geq 1$ ($e\geq 2$ and $n\geq 2$).
Let us thus consider two generic members $H_1$ and $H_2$ of such a pencil. By induction hypothesis, the varieties $X\cap H_i$, for $i=1,2$, 
are irreducible and satisfy the induction hypothesis. Their projection from $l$ is therefore $(n-1)$-dimensional.
The space ${\P}^{r}$ is spanned by $<H_1,H_2>$ and $X$ is non degenerate, thus $X^{\prime\prime}$ is $n$-dimensional.
\end{proof}

\begin{theorem} \label{th:3}Let $S$ be a non-degenerate (possibly singular) surface of degree $d$ and codimension $e\geq 2$ in ${\P}^r$.
Suppose that $S$ has a extremal secant line $l$, meeting $S$ along smooth points. Suppose that $S$ is the regular projection of a smooth
rational normal scroll $\overline{S}$ of degree $d$ in ${\P}^{d+1}$. Then, $S$ is a smooth rational scroll. 
\end{theorem}

\begin{proof}

Let $\Lambda$ denote the $(d-r)$-plane center of projection $\pi_{\Lambda}:{\P}^{d+1}\dashrightarrow {\P}^{r}$ mapping $\overline{S}$ to $S$.
By assumption, we have $\Lambda \cap\overline{S}=\emptyset$. Thus, no curves on $\overline{S}$ can be contracted to a point by $\pi_{\Lambda}$.
The $(d-r+2)$-plane $\overline{\Lambda}:=<\Lambda ,l>$ is thus a $k$-secant plane to $\overline{S}$, for some integer $k\geq 1$. 
\begin{lemma} The $(d-r+2)$-plane $\overline{\Lambda}$ is a $(d-r+3)$-secant plane.
Moreover, we have 
\[
\mu_{x}(\Lambda\cap \overline{S})\leq 2 \quad\forall x\in \overline{S}
\]
\end{lemma}
\begin{proof}
Let $\{\overline{p}_1 ,\cdots , \overline{p}_m \}$ denote the support of $\Lambda\cap \overline{S}$. 
Let $C$ be a hyperplane section of $S$ by a generic hyperplane of ${\P}^{r}$ containing $l$. 
Since $l$ meets $S$ at smooth points, $C$ is irreducible and $l$ is an extremal secant line to $C$.
By theorem \ref{th:GLP}, $C$ is a smooth rational curve of degree $d$. Let $\overline{C}:=<\Lambda , C>\cap \overline{S}$;
it is a rational normal curve of degree $d$, projecting  onto $C$ by $\pi_{\Lambda}$.
By construction, the $0$-scheme $\overline{\Lambda}\cap \overline{S}$ is a subscheme of $\overline{\Lambda}\cap \overline{C}$.
We have $\mu_{\overline{p}_i}(\Lambda\cap \overline{S})\leq \mu_{\overline{p}_i}(\Lambda\cap \overline{C})\leq 2$. 
Indeed, suppose to the contrary that $\mu_{\overline{p}_i}(\Lambda\cap \overline{C}) >2.$
Then, $\Lambda$ contains the projective tangent line $\mathbb{T}_{\overline{p_i}}(\overline{C})$ to $\overline{C}$ at $\overline{p_i}$,
so that we would get 
\[
\mu_{\overline{p}_i}(\Lambda\cap \overline{C})=\mu_{\overline{p}_i}(\mathbb{T}_{\overline{p_i}}(\overline{C})\cap \overline{C})\leq 2,
\]
for $\overline{C}$ has no $3$-secant line.

\begin{claim} Let $p_i:=\pi_{\Lambda}(\overline{p_i})$.
We have
\[
\mu_{\overline{p_i}}(\overline{\Lambda}\cap \overline{C})\geq \mu_{p_i}(C\cap l)
\]

\end{claim}
\begin{proof}
The claim is trivially true if $\mu_{p_i}(C\cap l)=1$.
Suppose that $\mu_{p_i}(C\cap l)\geq 2$.  The projectivized tangent line $\mathbb{T}_{\overline{p_i}}(\overline{C})$ is the
 limiting position of the lines $<\overline{p_i} ,\overline{x}>$, for $x\in \overline{C}\setminus \{ x \}$. Such a line cannot be contracted to a 
point by $\pi_{\Lambda}$, since $p_i$ is a smooth point of $C$.  
Therefore, the limiting position of $\pi_{\Lambda}( <\overline{p_i} ,\overline{x}>)$ is $\mathbb{T}_{p_i}(C)$.
By assumption, we have $l=\pi_{\Lambda}(\mathbb{T}_{\overline{p_i}}(\overline{C}))$, so that $\mu_{\overline{p_i}}(\overline{\Lambda}\cap \overline{C})=2$.
Let us show that we always have $\mu_{p_i}(C\cap l)\leq 2$.
Up to a change of projective coordinates, we may assume that $p_i:=(0:0:\cdots :1)$ in $<C>=Proj(\C [x_0,\cdots ,x_{r-1}])$ and that 
the saturated ideal defining $l$ in $<C>$ is $(x_0,\cdots ,x_{r-3})$. We write $<\overline{C}>=Proj(\C [x_{r},\cdots,x_{d},x_0,\cdots ,x_{r-1}])$ and 
denote by $I_{\overline{C}}$ the saturated ideal defining $\overline C$ in $<\overline{C}>$.
Let $\{ \overline{f}_1, \cdots ,\overline{f}_u\}$ be a Gr\" obner basis of $I_{\overline{C}}$ for lexicographic order on  
$\C [x_{r},\cdots,x_{d},x_0,\cdots ,x_{r-1}]$.
Since $\overline{C}$ is $2$-regular, all the $f_i$'s are degree $2$-polynomials. Let us now proceed to the elimination; the subset 
$\{ \overline{f}_{i_1}, \cdots ,\overline{f}_{i_v}\}$ of $\{ \overline{f}_1, \cdots ,\overline{f}_u\}$ consisting of quadric polynomials in the variables
$x_0, \cdots ,x_{r-1}$ only is a Gr\" obner basis of $I_C$ for the lexicographic order on $\C [x_0,\cdots ,x_{r-1}]$.
For all $j\in\{ 1,\cdots ,v\}$, let us write $\overline{f}_{i_j}=f^{0}_{i_j}+f^{1}_{i_j}$, where $f^{0}_{i_j}$ is a homogeneous polynomial in the variables
$x_{r-2}$ and $x_{r-1}$ only. Then $(f^{0}_{i_1},\cdots ,f^{0}_{i_v})$ is a defining ideal for the $0$-scheme $C\cap l$ in $l=Proj(\C [x_{r-2},x_{r-1}])$; 
this scheme is supported at the point $(0:1)$, hence the length of this $0$-scheme satisfies
\[
\mu_{p_i}(l\cap C)=min\{ deg_{x_{r-2}} ( f^{0}_{i_j})| j\in\{ 1 ,\cdots , v\}\,\,\text{and} \,\, f^{0}_{i_j}\not =0\},
\]
where  $deg_{x_{r-2}}(f^{0}_{i_j})$ denotes the degree of the polynomial  $f^{0}_{i_j}$, viewed as a polynomial in the variable $x_{r-2}$.
Since $deg_{x_{r-2}}(f^{0}_{i_j})\leq deg(f^{0}_{i_j})=2,\,\,\forall j\in\{ 1,\cdots , v\}$, we find $\mu_{p_i}(l\cap C)\leq 2$. 
This finishes the proof the claim.
\end{proof} 
Summing up the inequalities of the claim over the points of the support of the $0$-scheme $\overline{\Lambda}\cap \overline{C}$,
we find $k\geq length(\Lambda\cap \overline{C})\geq d-r+3$. 
From theorem \ref{th:sec}, we have $k\leq d-(d+1-2)+d-r+2=d-r+3$.
\end{proof}
Let $s$ denote the number of points $y$ of $\overline{\Lambda}:=<\Lambda,l>$, such that $\mu_{y}(\overline{\Lambda}\cap\overline{S})=2$.
We have $0\leq s\leq \frac{d-r+3}{2}$. Let $\{\overline{p}_1,\cdots \overline{p}_{\delta-s}\}$ denote the support of the $0$-scheme 
$\overline{\Lambda}\cap\overline{S}$. We may assume that $\mu_{\overline{p}_i}(\overline{\Lambda}\cap\overline{S})=2$ for $i\in\{1,\cdots ,s\}$
and  $\mu_{\overline{p}_i}(\overline{\Lambda}\cap\overline{S})=1$ otherwise.

\begin{lemma}\label{lem:1} Let $S$ be a surface satisfying the assumptions of theorem \ref{th:3}. Suppose that $\delta-s \not=3$ and $(\delta ,s)\not= (4,2)$. 
Then, $S$ is smooth.
\end{lemma}
\begin{proof}
Let $x$ be a point of $S$. Let $\phi$ denote the regular map $\overline{S}\xrightarrow{\pi_{\Lambda}|_{\overline{S}}}S$. 
If $x\in S$ does not belong neither to the double locus  nor to the ramification locus of $\phi$, $x$ is a smooth point of $X$. 
Suppose that $S$ is singular. Let $x$ be a singular point of $S$. 
Suppose first that the support of the fiber of $\phi$ at $x$ contains two distinct points $x_1$ and $x_2$.
Let $\overline{C}$ be a section of $\overline{S}$ by some hyperplane, generic among hyperplanes containing $\overline{\Lambda}:=<\Lambda,l>$.
We have $ \overline{\Lambda}:=<\mathbb{T}_{\overline{p}_{1}}(\overline{C}),\cdots ,\mathbb{T}_{\overline{p}_{s}}(\overline{C}), 
\overline{p_{\delta-2s+1}},\cdots ,\overline{p_{\delta-s}}>$.
Consider the following linear space 
\[
L_x:=\begin{cases}
<\mathbb{T}_{\overline{p}_{1}}(\overline{C}),\cdots ,\mathbb{T}_{\overline{p}_{s}}(\overline{C}), \overline{p_{\delta-2s+1}},\cdots ,\overline{p_{\delta-s-2}},
\overline{x_1},\overline{x_2}>\,\text{if}\, s\not= \delta /2\\
<\mathbb{T}_{\overline{p}_{1}}(\overline{C}),\cdots ,\mathbb{T}_{\overline{p}_{s-1}}(\overline{C}), \overline{x_1},\overline{x_2}> \,\text{if}\, s= \delta /2
\end{cases}
\]
The linear space $\pi_{\Lambda}(L_x)\subset <l,x>$, is a line meeting $l$ at a single point, since $x\not\in l$. 
Therefore, $\delta-s=3$ or $\delta=4$ and $s=2$. 
By assumption, this is not the case, hence the fiber of $\phi$ at a singular point of $S$
is supported at single point $\overline{x}$. 

Consider a section $\overline{D_x}$ of $\overline{S}$ by some hyperplane $\overline{H}_x$ in ${\P}^{d+1}$, 
generic among the hyperplanes containing $<\Lambda,l>$. Since $\overline{H}_x$ contains 
$\Lambda$, $H_x:=\pi_{\Lambda}(\overline{H}_x)$ is a hyperplane in ${\P}^r$. 
By Bertini's smoothness theorem, since $d+1-(d-r+4)=e-1\geq 1$, the degree $d$ curve $\overline{D_x}$ is smooth hence irreducible.
The morphism $\phi_{\overline{D_x}}$ maps $\overline{D_x}$ to a curve $D_x$ on $S$, 
hyperplane section of $S$ be some hyperplane containing $<l,x>$; hence, $deg(D_x)=d$.
The curve $D_x$ is the linear projection of an irreducible curve, hence is irreducible.

If $x$ doesn't belong to the ramification locus of $\phi|_{\overline{D_x}}$, $x$ is a smooth point of $D_x$.
Since $D_x$ is a Cartier divisor of $S$, $x$ is a smooth point of $S$.
Therefore, $x$ is a ramification point of $\phi|_{\overline{D_x}}$, i.e. $\mathbb{T}_{\overline{x}}(\overline{D_x})$ meets $\Lambda$.

Consider the following linear space
\[
L_x:=\begin{cases}
<\mathbb{T}_{\overline{p}_{1}}(\overline{D_x}),\cdots ,\mathbb{T}_{\overline{p}_{s}}(\overline{D_x}), 
\overline{p_{\delta-2s+1}},\cdots ,\overline{p_{\delta-s-2}}, \mathbb{T}_{\overline{x}}(\overline{D}_x)>\, ,\text{if}\quad  s\not= \frac{\delta}{2} \\
<\mathbb{T}_{\overline{p}_{1}}(\overline{D_x}),\cdots ,\mathbb{T}_{\overline{p}_{s-1}}(\overline{D_x}),\mathbb{T}_{\overline{x}}(\overline{D}_x)>,
\, \text{if}\quad  s= \frac{\delta}{2}  
\end{cases}
\] 
As in the previous case, the linear space $\pi_{\Lambda}(L_x)\subset <l,x>$, is a line meeting $l$ at a single point. 
Therefore, $\delta-s=3$ or $\delta=4$ and $s=2$. This contradicts our assumption. Thus, $S$ is smooth.
\end{proof}

\begin{lemma} \label{lem:2} Let $S$ be a surface satisfying the assumption of theorem \ref{th:3}. Suppose that $S$ is singular and $e\geq 3$.
\begin{enumerate}
\item Then, $S^{\prime\prime}$, its projection from $l$, is a surface $S(0,r-3)$, cone over a rational normal curve of degree $r-3$.
All the singularities of $S$ lie in the $2$-plane projecting from $l$ onto the vertex of the cone $S^{\prime\prime}$. In particular,
$S$ has isolated singularities.
\item Let $s$ denote the number of points $x\in \overline{\Lambda}\cap \overline{S}$ for which
$\mu_{x}(\overline{\Lambda}\cap \overline{S})=2$. Then, we have $s\in 2 \mathbb{Z}_{\geq 1}$ and $\overline{S}=S(s/2,d-s/2)$.
\end{enumerate}
\end{lemma}
\begin{proof}
Let us focus on the possible singularities of $S$.
If $x\in S$ doesn't belong neither to the double locus nor to the ramification locus of $\pi_{\Lambda}$, $x$ is a smooth point of $S$, 
since $\overline{S}$ is smooth.  

Suppose that two points $\alpha$ and $\beta$ get contracted to a double point $x$ of $S$. 
If the linear system $\mathcal{L}_{l,x}$ of hyperplane section of $S$ by  hyperplanes of ${\P}^{r}$ through $<l,x>$
satisfy the assumption of theorem \ref{th:irr}, a general member of $\mathcal{L}_{l,x}$ is a singular irreducible curve having an extremal secant line. 
This contradicts theorem \ref{th:GLP}. Therefore, the projection of $X$ from the $2$-plane  $<l,x>$ is a curve $C$ in ${\P}^{r-3}$, i.e. the surface $S^{\prime\prime}:=\pi_{l}(S)$ is the cone $<\pi_{l}(x),C>$ over a rational normal curve $C$ of degree $r-3$.

Suppose that $<\Lambda ,y>$ contains a tangent line $T$ to $\overline{S}$ at $y$, i.e $x :=\pi_{\Lambda}(y)$ belongs to the ramification 
locus of $\pi_{\Lambda}$. We write $\overline{\Lambda}:=<\Lambda ,l>$.
Then, unless the linear system of hyperplane sections through $<\overline{\Lambda},y>$ 
doesn't satisfy assumptions of theorem \ref{th:irr}, i.e. $X^{\prime\prime}=<x,C>$ for a rational normal curve $C$ of degree $r-3$,
a generic hyperplane section through $<\overline{\Lambda} ,y>$ is a rational normal curve $D$ tangent to $x$ along the line $T$.
Therefore, $x$ is a smooth point of the linear section of $X$, $\pi_{\Lambda}(D)$, for this curve has an extremal secant line $l$.    
Thus, $x$ is a smooth point of $S$.
In summary, if $x$ is a singular point of $S$, then $S^{\prime\prime}$ is a cone $<q, C>$ over a rational normal curve of degree $r-3$, and
$\pi_{l}(x)=q$.

Let $\{ \overline{p}_1 ,\cdots ,\overline{p}_{m}\}$ denote the support of $\overline{\Lambda}\cap X$, where $m\leq d-e+1$. 
If $s$ denotes the number of points $x\in\overline{\Lambda}\cap \overline{S}$ for which $\mu_{x}(\overline{\Lambda}\cap \overline{S})=2$,
we have $m=\delta -s$, where $\delta:=d-r+3$.
Since $\overline{S}$ is smooth, for $i\not =j \in \{ 1,\cdots ,m\}$, we have $\overline{p}_j\not\in \overline{f}_{i}$, 
where $\overline{f}_{i}$ denotes the fiber of $\overline{S}$ passing through $\overline{p}_{i}$
 
The projection $\pi_{\overline{\Lambda}}$ from $\overline{S}$ to $S^{\prime\prime}$, induces an elementary transformation at 
$\overline{p}_1 ,\cdots ,\overline{p}_{m}$ from $\overline{X}$ to the surface $\mathbb{F}_{r-3}$, blowing up of the cone $S^{\prime\prime}$ at its vertex $q$. 

Indeed, the indeterminacy locus of the rational map induced by $\pi_{\overline{\Lambda}}$ on $\overline{S}$
is resolved by blowing up $\overline{p}_1 ,\cdots , \overline{p}_{m}$. 
Let $\tilde{S}$ denote the blowing up of $X$ at $\overline{p}_1 ,\cdots ,
\overline{p}_{m}$ and, for $i=1,\cdots ,m$, let $E_i$ denote the exceptional curve over the point $p_i$ and 
$\tilde{f_i}$ the strict transform of the fiber of $\overline{f}_{i}$.

The resolving regular map $\phi$ from $\tilde{S}$ to $S^{\prime\prime}$ is associated to  a sub-system of the complete linear system $\|H-\sum_{i=1}^{m} E_i\|$,
so that it contracts $\tilde{f}_i$ to a point of ${X}^{\prime\prime}$ and sends $E_i$ to a line of $X^{\prime\prime}$.

Since $deg(C)\geq 2$, the only lines of $<q,C>$ are lines joining $q$ to a point of $C$, i.e. blow-down of fibers of the scroll $\mathbb{F}_{r-3}$.
The birational map from $\overline{S}$ to $\mathbb{F}_{r-3}$, induced by $\pi_{\overline{\Lambda}}$ factors 
though an elementary transformation $\sigma$ of $\overline{S}$ at $\overline{p_1} ,\cdots ,\overline{p_m}$ and a composition of 
blow-downs of $(-1)$-curves to $\mathbb{F}_{r-3}$.
Let  $\Sigma$ be the smooth rational ruled surface image of $\overline{S}$ by  $\sigma$. The surface $\Sigma$ is isomorphic to a surface of type $\mathbb{F}_{n}$
for some integer $n\geq 0$. If $n\not =1$, by minimality of $\mathbb{F}_{n}$ we get $n=r-3$. If $n=1$, then $\Sigma$ has a unique $(-1)$-curve, 
which gets blown-down to a point of $\mathbb{F}_{r-3}$; thus, $\mathbb{F}_{0}\simeq\mathbb{P}^2$ is an elementary transformation of $\mathbb{F}_{r-3}$
at a point $y$ of $\mathbb{F}_{r-3}$, so that $r-4=0$ or $r-2=0$. From the assumption $e\geq 3$, neither of the two cases can happen.
In conclusion, $\Sigma$ is isomorphic to $\mathbb{F}_{r-3}$ and $\phi$ is an elementary transformation of $\overline{X}$ at $\overline{p_1},\cdots ,
\overline{p_{m}}$.

\begin{claim} Let $\psi$ denote the birational map from $\overline{S}$ to $X^{\prime\prime}=<q,C>$, induced by $\pi_{\overline{\Lambda}}$;
${\psi}^{-1}$ is defined away from the vertex of the cone $q$.
\end{claim}

\begin{proof}
If $x\in S$ is a smooth point, we denote by $f_x$ the image by $\pi_{\Lambda}$ of the fiber of $\overline{S}$ at the preimage of $x$ by $\pi_{\Lambda}$.
Since $\Lambda \cap \overline{S}=\emptyset$, no fibers of $\overline{S}$ can be contracted to a point by $\pi_{\Lambda}$; $f_x$ is thus a line.
We call $f_x$ the fiber of $S$ at $x$.

Since $\psi$ induces an elementary transformation at $\overline{p_1} ,\cdots ,\overline{p_m}$, 
it is enough to show that for all $i=1,\cdots ,m$ the fibers of $S$ at $p_i$ are contracted by $\psi$ to $q$.
 
Suppose, to the contrary, that there exist $i\in \{ 1,\cdots ,m\}$ such that the fiber of $S$ at $p_{i}$
is contracted to $z\not =q$ on $<q,C>$. Pick a generic point $y$ on the line $<q,z>$;
since $y\not= q$, $y$ is the image by $\pi_{l}$ of a smooth point $x$ of $S$. Since $y$ is generic on $<q,z>$, for all $j\in \{1,\cdots ,m\}$ we have 
$f_x \not= f_{p_{j}}$, so that $f_x$ meets the $2$-plane $<l,f_{p_{i}}>$ transversally at some point $x^{\prime}$.
Since $\pi_{l}(x^{\prime})=z\not =q$, $x^{\prime}$ is a smooth point of $S$, hence cannot lie on the fiber $f_{p_{i}}$.
Remark that the generic member of the linear system of hyperplane sections through $l$ and $x^{\prime}$ is not irreducible.
Thus, $\pi_{z}(<q,C>)$ is a curve $C^{\prime}$ of $\mathbb{P}^{r-3}$. The cone $<q,C>$  thus coincides with the cone $<z,C^{\prime}>$.
Since $deg(C)=r-3\geq 2$, $q=z$ and we get a contradiction.
\end{proof}

Let $\phi$ be the birational inverse of $\psi$; it is an elementary transformation at $m$ distinct
points $q_1 ,\cdots ,q_{m}$ of $\mathbb{F}_{r-3}$. 

Since the inverse of $\psi$ is defined away from the vertex of cone $<q,C>$, the $m$ points $q_1 ,\cdots ,q_{m}$ 
must lie on the unique $(-(r-3))-$ curve  on $\mathbb{F}_{r-3}$, $E_q$, 
the exceptional locus of the  blow-up of the cone $<q,C>$ at its vertex $q$. 
By the blowing-down of the fibers $g_i$ of the points $q_i$ for 
$i=1,\cdots ,m$, $E_q$ is blown down to a $(-(r-3+m))-$ irreducible curve on $\overline{S}$.
We write $\overline{S}:=S(\alpha,\beta)$, with $1\leq \alpha\leq \beta$.
Since  $\overline{S}\simeq \mathbb{F}_{\beta -\alpha}$, we find $\beta-\alpha=r-3+m$.
From $\alpha +\beta =d$, we find $s\in 2\,\mathbb{Z}_{\geq 1}$ and $\alpha=s/2$
$\beta=d-s/2$. 

The singular locus of $S$ lies in the $1$-dimensional intersection of $S$ with the $2$-plane $<l,q>$;
since $l\cap S$ is supported at smooth points of $S$, $Sing(S)$ is $0$-dimensional.

\end{proof}
If  $\delta-s=3$, we find $s=0$ or $s=1$.
Thus, combining lemma \ref{lem:1} and \ref{lem:2}, we deduce that, if $S$ is singular, either $e\geq 3$, $\delta=4$ and $s=2$,
or $e=2$ and $\delta -s\in\{2,3\}$.
\begin{lemma} Suppose that $S$ satisfies the assumptions of theorem \ref{th:3}. If $e=2$, $S$ is smooth.
\end{lemma} 
\begin{proof} 
By lemma \ref{lem:1}, we can assume that $\delta-s=3$ or $\delta=4$ and $s=2$. That is to say
$(d,\delta ,s)\in\{ (5,4,2),(5,4,1),(4,3,0)\}$.
Suppose first that $(d,\delta ,s)=(4,3,0)$. Then, we have $r=4=d$.
By corollary \ref{cor:2}, $\pi_{l}(S)={\P}^2$. Hence for any point $x$ of $S$, $\pi_{l,x}(S)={\P}^1$.
In particular, the generic member $D_x$ of the linear system cut out on $S$ by hyperplanes through $<x,l>$ is irreducible, 
since generic hyperplane sections of ${\P}^1$ are irreducible. 
Suppose first that $(d,\delta ,s)=(4,3,0)$. Then, we have $r=4=d$, so that $D_x$ is a rational quartic in ${\P}^3$ having a $3$-secant line $l$.
By theorem \ref{th:GLP}, $x$ is a smooth point of $D_x$.
The point $x$ is thus a smooth point of $S$, for $D_x$ is a Cartier divisor on $S$.
The surface $S$ is thus smooth.

If $(d,\delta ,s)\in\{(5,4,2),(5,4,1)\}$, we find $d=5=r-1$, so that $D_x$ is either a possibly singular curve of degree $5$ in ${\P}^3$.
If $dim<D_x>=2$, $D_x$ is a plane quintic curve, so that any line meets $D_x$ with multiplicity $deg(D_x)=5> \delta =4$.
This gives a contradiction, since, by construction, $D_x\cap l=S\cap l$ as $0$-schemes.
Therefore, $dim<D_x>=3$, and  $D_x$ has an extremal secant line $l$.
The curve $D_x$ is thus smooth, showing that $x$ is a smooth point of $S$.
\end{proof}
The only remaining case is $e\geq 3$ and $(\delta ,s)=(4,2)$, i.e. $(d,\delta ,s)=(e+3,4,2)$ and $d\geq 6$.
Suppose that $S$ is singular. In particular, $\pi_{l}(S)=S(0,r-3)$.
We have $\overline{S}=S(1,d-1)$. So that there is a line $l^{\prime}$ on $\overline{S}$ meeting 
every fiber $f_x$ of $\overline{S}$ at a single point $y(x)$; re-parameterizing $\overline{S}$ by $l^{\prime}$, we can assume $x=y(x)$. 
Since $\Lambda$ doesn't contract any curve on $\overline{S}$,
$\pi_{\Lambda}(l^{\prime})$ is a line on $S$, meeting $\pi_{\Lambda}(f_x)$ for all $x\in l^{\prime}$.
Suppose that $dim(\pi_{l}(\pi_{\Lambda}(l^{\prime})))=1$, then $\pi_{l}(\pi_{\Lambda}(l^{\prime})))$ is a line on $S(0,r-3)$ meeting
all lines  $\pi_{\overline{\Lambda}}(f_x)$ for all $x\in l^{\prime}\setminus l$. 
Since $r-3\geq 2$,i.e. $S(0,r-3)\not={\P}^2$; this is impossible.
Therefore, $l^{\prime}$ gets contracted by $\pi_{\overline{\Lambda}}$ to  a point on the cone $S(0,r-3)$.
Let us consider the elementary transformation at $\{\overline{p}_1,\overline{p}_2 \}$ from $S(1,d-1)$ to ${\mathbb{F}}_{r-3}$ induced by 
$\pi_{\overline{\Lambda}}|_{\overline{S}}$.
Let $\tilde{S}$ denote the blow up of $\overline{S}$ at $\{\overline{p}_1,\overline{p}_2 \}$ and $E_1,E_2$ the exceptional divisors on $\tilde{S}$ over
$\overline{p}_1$ and $\overline{p}_2$ respectively.
We have $Pic(\tilde{S})=\mathbb{Z}H^{*}\oplus \mathbb{Z}f^{*}\oplus \mathbb{Z}E_1\oplus \mathbb{Z}E_2$, where $H^{*}$ is the pull-back of the class
$H$ of hyperplane sections in $S(1,d-1)$ and $f^{*}$ is the pull-back of the class of a fiber in $S(1,d-1)$.
Recall that $[l^{\prime}]=H-(d-1)f$ in $Pic(\overline{S})$, so that the class of the strict transform $\tilde{l^{\prime}}$ of $l^{\prime}$ in $\tilde{S}$
is of the form $H^{*}-(d-1)f^{*}-\epsilon_{1} E_1 -\epsilon_{2}E_2$, with $\epsilon_{i}\in\{0,1\}$ for $i=1,2$.
Moreover, since $\pi_{\Lambda}(l^{\prime})\not =l$, we get $\epsilon_1 +\epsilon_{2}\leq 1$.
If $\psi$ contracts $l^{\prime}$ to a point, we have $(\tilde{l^{\prime}})^2=-1$. Since $(\tilde{l^{\prime}})^2=-d+2+\epsilon_1 +\epsilon_{2}$,
this implies that $d\in\{3,4\}$. Since $d=e+3\geq 6$, $\psi$ cannot contract $l^{\prime}$. Therefore, since $\pi_{\overline{\Lambda}}(l^{\prime})$
is a point on $S(0,r-3)$, $\psi (l^{\prime})= \gamma$, where $\gamma$ is the unique irreducible curve on $\mathbb{F}_{r-3}$,
such that $\gamma ^{2}=-(r-3)$.
Therefore $\tilde{\gamma}=\tilde{l^{\prime}}$. We find, $-(r-3)+2=-d+1+\epsilon_1 +\epsilon_2$, so that 
 \[
\epsilon_1 +\epsilon_2 -1=d-(r-3)=\delta =4.
\]
This contradicts trivially $\epsilon_1 +\epsilon_{2}\leq 1$. Therefore, $S$ is smooth.
\end{proof}

\begin{corollary} \label{cor:class} Let $X$ be a non-degenerate complex complex projective variety of degree $d$ and codimension $e\geq 2$ in ${\P}^r$.
Suppose that $X$ has an extremal secant line $l$. Then $X$ is either
\begin{itemize}

\item a cone $<L,V>$, where $L$ is a linear space of dimension $k\geq -1$ and $V$ is the Veronese surface in ${\P}^5$,

\item a cone $<L,V^{\prime}>$, where $L$ is a linear space of dimension $k\geq -1$ and $V^{\prime}$ is the isomorphic projection of the Veronese surface to 
${\P}^4$,
\item a cone $<L,X_0>$, where $L$ is a linear space of dimension $k\geq -1$ and $X_0$ is either a smooth rational scroll with an extremal secant line
or a smooth rational curve with an extremal secant line,

\end{itemize}

\end{corollary}

\begin{proof}
If $X$ is smooth this is a consequence of our previous results \cite{B}. 

Suppose that $X$ is singular. 

By theorem \ref{th:2}, we can  assume that there exists an extremal secant line $l$ to $X$ meeting $X$ along smooth points.
Applying theorem \ref{th:B00}, we deduce that $X$ is the regular projection of a variety of minimal degree $\overline{X}$ in ${\P}^{d+n-1}$.

Let us first assume that $\overline{X}$ is smooth. By theorem \ref{th:min} ($e\geq 2$), $\overline{X}$ 
is either the Veronese surface $V$ in ${\P}^5$ or a smooth rational normal scroll.

\begin{lemma}\label{lem:class1} Let $X$ satisfying the assumption of theorem \ref{cor:class}. Suppose that $\overline{X}$ is smooth. 
Then $X$ is either the Veronese variety $V$ in ${\P}^5$, $V^{\prime}$ its isomorphic projection to ${\P}^4$, 
or a smooth rational scroll $X_0$ having an extremal secant line. 
\end{lemma}

\begin{proof}
Let us first assume that $\overline{X}=V$. Suppose $\overline{X}\not = X$. The regular projection to $X$ in ${\P}^{r}$, is either an isomorphism 
and $X=V^{\prime}$ or $X$ is the Steiner surface in ${\P}^{3}$. The last case can be ruled out, since $e\geq 2$.

Suppose now that $\overline{X}$ is a smooth rational normal scroll and assume that $X\not = \overline{X}$.
If $X$ is singular, the singularities of $X$ were created by the projection. 

The lemma is true if $n=2$ (theorem \ref{th:3}), so we can assume $n\geq 3$. 
Let us prove by induction on $n\geq 2$ that $X$ is smooth. 
Suppose that the lemma is true for $n-1\geq 2$. 
Let $x$ be a point of $X\setminus l$. We wish to show that $x$ is a smooth point of $X$. Let $h_1$ be an hyperplane section of $X$ 
by a generic hyperplane $H$ through $<l,x>$.

Suppose that $h_1$ is not irreducible; then, $\pi_{\pi_{l}(x)}(X^{\prime\prime})$ is a curve $C$, so that $X^{\prime\prime}$ is the 
$2$-dimensional cone $<\pi_{l}(x), C>$.
This contradicts corollary \ref{cor:2}; the hyperplane section $h_1$ is thus irreducible. 

The line $l$ is an extremal secant line to $h_1$. Moreover, the variety $h_1$ is the projection of the smooth rational normal 
scroll $\overline{h_1}$, hyperplane section of $\overline{X}$ by the hyperplane $<\Lambda ,H>$. 
By induction hypothesis, $h_1$ is smooth, so that $x$ is a smooth point of $h_1$.
Since $h_1$ is a Cartier divisor on $X$, the point $x$ is a smooth point of $X$.  
\end{proof}

Only remains the case $\overline{X}$ singular. By theorem \ref{th:min}, $\overline{X}$ is a cone $<L,\tilde{X}>$, 
where $L$ is a linear space of arbitrary dimension and $\tilde{X}$ is a smooth rational normal scroll. 
Consider again the center of projection $\Lambda$, of the regular projection mapping $\overline{X}$ to $X$. Notice that, since  $\Lambda \cap X=\emptyset$,
no lines $<p,q >$ joining a point of $p$ of $\tilde{X}$ to a point $q$ of the vertex of the cone $L$ can be contracted to a point by the projection to $X$. 

Thus, the cone $<L,\tilde{X}>$, projects by $\pi_{\Lambda}$ onto a cone $<L^{\prime},X_{0}>$, where the rational scroll $X_0$ 
is the image of $\tilde{X}$ by $\pi_{\Lambda}$. From the inequality
\[
dim(L^{\prime})+dim( X_0 )+1=dim(X)\leq dim(L)+dim(\tilde{X})+1=dim(\overline{X})=dim(X) 
\]
we deduce that $dim(L^{\prime})=dim(L)$. By proposition \ref{lem:class1}, if $X_0$ has an extremal secant line, 
$X_0$ is smooth since it is the projection of a smooth rational normal scroll. 

So we may assume that $l$ is not an extremal secant line to $X_0$, i.e. $l$ does not lie in the linear space $<X_0>$.
Since $l$ meets $X$ along smooth points of $X$, $l\cap L=\emptyset$. Moreover, we have $l\cap X_0=\emptyset$.

Indeed, suppose, to the contrary, that there exists a point $p \in l\cap X_0$. Since $l\not \subset <X_0>$, the line $l$ meets $X$ at another point $q$ of 
$X\setminus(L\cup X_0)$. The line $l=<q,p>$ therefore lies   in the $(k+1)$-plane $<p,L>$ and meets $L$ at some point $x\in L$. Thus, $l$ coincides with 
the ruling line $<p,x>\subset X$. This gives a contradiction.  

Therefore, our extremal secant line $l$ to $X$ is the projection of a curve in $<\Lambda ,l>$ meeting $\overline{X}$ away from $L$ and $\tilde{X}$.
Let us rule out this case. First, notice that from corollary \ref{cor:2} the variety $X$ projects from $l$ onto a variety $X^{\prime\prime}$ 
of minimal degree $d$ in ${\P}^{r-2}$. 

Let us first assume that $l\cap <X_0>\not = \emptyset$; we denote by $x$ the intersection point $l\cap <X_0>$. Notice that $\pi_l (X_0)$ is a 
variety of minimal degree. Indeed, the variety $\pi_l (X)$ is the join variety $<\pi_l (L),\pi_l (X_0 )>$, 
which turns out to be a variety of minimal degree of dimension $n$ and degree $d$ in ${\P}^r$. Since $l\cap L=\emptyset$, we can realize  $\pi_l$
as a projection onto the $(r-2)$-plane $<L,\pi_l (<X_0>)>$, so that $\pi_l (L)=L$ and $\pi_l (X)$ is the cone $<L,\pi_l (X_0)>$.
Therefore, $\pi_l (X_0)$ is a variety of degree $d$, dimension $dim(X_0)=n-dim(L)$ and the same codimension in $<\pi_l (X_0)>$ as $\pi_l (X)$
in $<\pi_l (X)>$, hence $\pi_l (X_0)$ is a variety of minimal degree. 
Let $q$ be a generic point on $l$, we can realize $\pi_q$ as a projection from ${\P}^{r}$ onto the $(r-1)$-plane $<<X_0>, \pi_q (L)>$.
Therefore, its restriction to $X_0$ is the identity map. Since $\pi_l =\pi_x \circ \pi_q $,  the projection $\pi_x$ 
induces a birational map from $X_0$ onto its image $\pi_l (X_0)$. Since $\pi_l (X_0)$ is a variety of minimal degree, it is linearly normal, 
so that $X_0 =\pi_l (X_0)$, i.e. $x\not \in <X_0>$. 
 
Since $l\cap <X_{0}>=\emptyset$, the the projection $\pi_{l}$ induces the identity map between $X_0$ and its image $X_{0}^{\prime\prime}$, 
which is a variety of minimal degree $d$. Therefore, $X$ is a variety of variety of minimal degree. 
This contradicts $\overline{X}\not = X$.

\end{proof}

Since a cone has the same regularity as its base and smooth scrolls satisfy the regularity conjecture \cite{B}, if $e\geq 2$, 
varieties having an extremal secant line satisfy the regularity conjecture.

\textbf{Acknowledgements} I am grateful to Atsushi Noma for pointing out this mistake in my previous paper.

\end{document}